\documentclass[11pt]{article}
\PassOptionsToPackage{obeyspaces}{url}

\usepackage{amsfonts,amsmath,amssymb,amsthm}
\usepackage{mathtools}
\usepackage{enumitem}
\usepackage{graphicx}
\usepackage{tikz}
\usetikzlibrary{arrows.meta,calc}
\usepackage{bm}
\usepackage{cancel}
\usepackage{comment}
\usepackage{etoolbox}
\usepackage{xstring}

\usepackage[hidelinks]{hyperref}
\usepackage{bookmark}
\bookmarksetup{open,numbered,depth=5}
\usepackage[nameinlink]{cleveref}

\allowdisplaybreaks[4]

\patchcmd{\thebibliography}
{\leftmargin\labelwidth}
{\leftmargin\labelwidth\addtolength\itemsep{-0.3\baselineskip}}
{}{}

\theoremstyle{plain}
\newtheorem{theorem}{Theorem}[section]

\newtheorem{lemma}[theorem]{Lemma}

\newtheorem{conjecture}[theorem]{Conjecture}

\theoremstyle{definition}

\newtheorem{problem}[theorem]{Problem}

\crefname{fact}{fact}{facts}
\Crefname{fact}{Fact}{Facts}
\crefname{claim}{claim}{claims}
\Crefname{claim}{Claim}{Claims}
\crefname{observation}{observation}{observations}
\Crefname{observation}{Observation}{Observations}
\crefname{question}{question}{questions}
\Crefname{question}{Question}{Questions}
\crefname{construction}{construction}{constructions}
\Crefname{construction}{Construction}{Constructions}
\crefname{exercise}{exercise}{exercises}
\Crefname{exercise}{Exercise}{Exercises}
\crefname{problem}{problem}{problems}
\Crefname{problem}{Problem}{Problems}

\DeclareMathOperator{\ex}{ex}

\newcommand*{\arXiv}[1]{\href{https://arxiv.org/abs/#1}{arXiv:#1}}

\definecolor{mred}{RGB}{217,38,38}
\definecolor{mgreen}{RGB}{38,217,38}
\definecolor{mblue}{RGB}{38,38,217}
\definecolor{mcyan}{RGB}{38,128,128}
\definecolor{mmagenta}{RGB}{128,38,128}
\definecolor{myellow}{RGB}{128,128,38}
\definecolor{mgray}{RGB}{98,98,98}

\title{Exact Tur\'an densities of even zycles}
\author{Shuaichao Wang\thanks{\raggedright%
		Faculty of Mathematics and Statistics, Central China Normal University,
		Wuhan 430079, PR China. \texttt{shuaichao@ccnu.edu.cn}.}}
\date{}
\hypersetup{pdftitle={Exact Turan densities of even zycles},pdfauthor={Shuaichao Wang}}

\begin{document}
	\maketitle
	
	\begin{abstract}
		
		We determine the exact Tur\'an density
		of the $3$-uniform zycle $Z_\ell$ for every
		even integer $\ell\ge 4$ with $3\nmid\ell$,
		showing that $\pi(Z_\ell)=5/9$.
		More generally, we prove that
		$\pi(Z_\ell)\le 5/9$ for all even $\ell\ge 4$.
		In the range of equality, we also establish
		a stability theorem, showing that every almost extremal
		$Z_\ell$-free $3$-graph is close in edit distance
		to the classical Tur\'an construction.
		
	\end{abstract}
	
	\section{Introduction}\label{sec:introduction}
	
	For an integer $r\ge2$, an \emph{$r$-uniform hypergraph}, or
	\emph{$r$-graph}, is a pair $H=(V(H),E(H))$ with
	$E(H)\subseteq\binom{V(H)}r$. Write $e(H)=|E(H)|$.
	Given a family $\mathcal F$ of $r$-graphs, an $r$-graph is
	\emph{$\mathcal F$-free} if it contains no member of
	$\mathcal F$ as a subgraph. The \emph{Tur\'an number}
	$\ex(n,\mathcal F)$ is the maximum number of edges in an
	$\mathcal F$-free $r$-graph on $n$ vertices, and the
	\emph{Tur\'an density} of $\mathcal F$ is
	\[
	\pi(\mathcal F)=\lim_{n\to\infty}
	\frac{\ex(n,\mathcal F)}{\binom nr}.
	\]
	The limit exists because averaging over vertex subsets shows
	that the normalized extremal numbers are nonincreasing.
	For $\mathcal F=\{F\}$, we write $\ex(n,F)$ and $\pi(F)$.
	All subgraphs in this paper are ordinary, not necessarily
	induced, subgraphs.
	
	Determining Tur\'an densities is a central problem in extremal
	combinatorics; see Keevash's survey~\cite{KeevashSurvey}
	for an overview. For hypergraphs of uniformity at least three,
	even small forbidden configurations can be difficult to handle.
	The classical tetrahedron problem asks whether
	$\pi(K_4^{(3)})=5/9$, where $K_4^{(3)}$ is the complete
	$3$-graph on four vertices. Despite substantial progress,
	exact densities are known for relatively few individual
	$3$-graphs. Classical examples include the Fano plane
	$\mathbb F$, the $3$-graph on seven vertices in which every
	pair belongs to exactly one edge, and the $3$-graph $F_{3,2}$
	on $[5]$ with edge set $\{123,145,245,345\}$, where $ijk$
	abbreviates $\{i,j,k\}$. De Caen and F\"uredi~\cite{dCF}
	proved that $\pi(\mathbb F)=3/4$, while F\"uredi,
	Pikhurko and Simonovits~\cite{FPS} established
	$\pi(F_{3,2})=4/9$.
	Another classical example is the $3$-graph $F_5$ with edge
	set $\{123,124,345\}$, for which Frankl and
	F\"uredi~\cite{FF} proved that $\pi(F_5)=2/9$.
	For $3$-uniform expansions of complete graphs of order
	at least four, Mubayi~\cite{MubayiExpansion} determined
	the Tur\'an densities, and Pikhurko~\cite{PikhurkoExpansion}
	subsequently determined the extremal numbers for all
	sufficiently large orders.
	
	More recently, exact results have been obtained for tight
	cycles. For $m\ge4$, let $C_m^{(3)}$ be the $3$-graph on
	$\mathbb Z_m$ with edges $\{i,i+1,i+2\}$ for
	$i\in\mathbb Z_m$. Kam\v{c}ev, Letzter and
	Pokrovskiy~\cite{KLP} determined their densities for all
	sufficiently large lengths not divisible by three.
	Bodn\'ar, Le\'on, Liu and Pikhurko~\cite{BLLP}
	extended this result by proving that
	$\pi(C_m^{(3)})=2\sqrt3-3$ for every $m\ge7$ with
	$3\nmid m$.
	
\begin{figure}[htbp]
		\centering
		\begin{tikzpicture}[
			font=\scriptsize,
			vertex/.style={circle,fill=black,inner sep=1.25pt},
			firstedge/.style={draw=teal!80!black,fill=teal!60,
				fill opacity=.16,line width=.65pt,line join=round},
			secondedge/.style={draw=violet!75!black,fill=violet!50,
				fill opacity=.11,line width=.65pt,line join=round,
				dash pattern=on 2.2pt off 1.3pt}]
			% Each regular polygon has the same 3.2 cm width and height.
			% Consecutive vertices x_i,y_i form the ith designated pair.
			\foreach \ell/\shift in {2/0,4/4.8,8/9.6} {
				\begin{scope}[xshift=\shift cm]
					\pgfmathtruncatemacro{\last}{\ell-1}
					\pgfmathsetmacro{\radius}{1.60/cos(90/\ell)}
					\pgfmathsetmacro{\labelradius}{1.85/cos(90/\ell)}
					\foreach \i in {0,...,\last} {
						\coordinate (x\i) at ({90+90/\ell-360*\i/\ell}:\radius);
						\coordinate (y\i) at ({90-90/\ell-360*\i/\ell}:\radius);
					}
					\foreach \i in {0,...,\last} {
						\pgfmathtruncatemacro{\j}{mod(\i+1,\ell)}
						\path[firstedge] (x\i)--(y\i)--(x\j)--cycle;
						\path[secondedge] (x\i)--(y\i)--(y\j)--cycle;
					}
					\foreach \i in {0,...,\last} {
						\node[vertex] at (x\i) {};
						\node[vertex] at (y\i) {};
						\node at ({90+90/\ell-360*\i/\ell}:\labelradius) {$x_\i$};
						\node at ({90-90/\ell-360*\i/\ell}:\labelradius) {$y_\i$};
					}
					\ifnum\ell=2
					\node[font=\small] at (0,-2.30) {$Z_2\cong K_4^{(3)}$};
					\else
					\node[font=\small] at (0,-2.30) {$Z_\ell$};
					\fi
				\end{scope}
			}
		\end{tikzpicture}
		\caption{The zycles $Z_2$, $Z_4$ and $Z_8$.
			Solid teal triangles represent $\{x_i,y_i,x_{i+1}\}$ and dashed
			purple triangles represent $\{x_i,y_i,y_{i+1}\}$, with indices
			modulo $\ell$. Only marked points are vertices.}
		\label{fig:zycles}
	\end{figure}
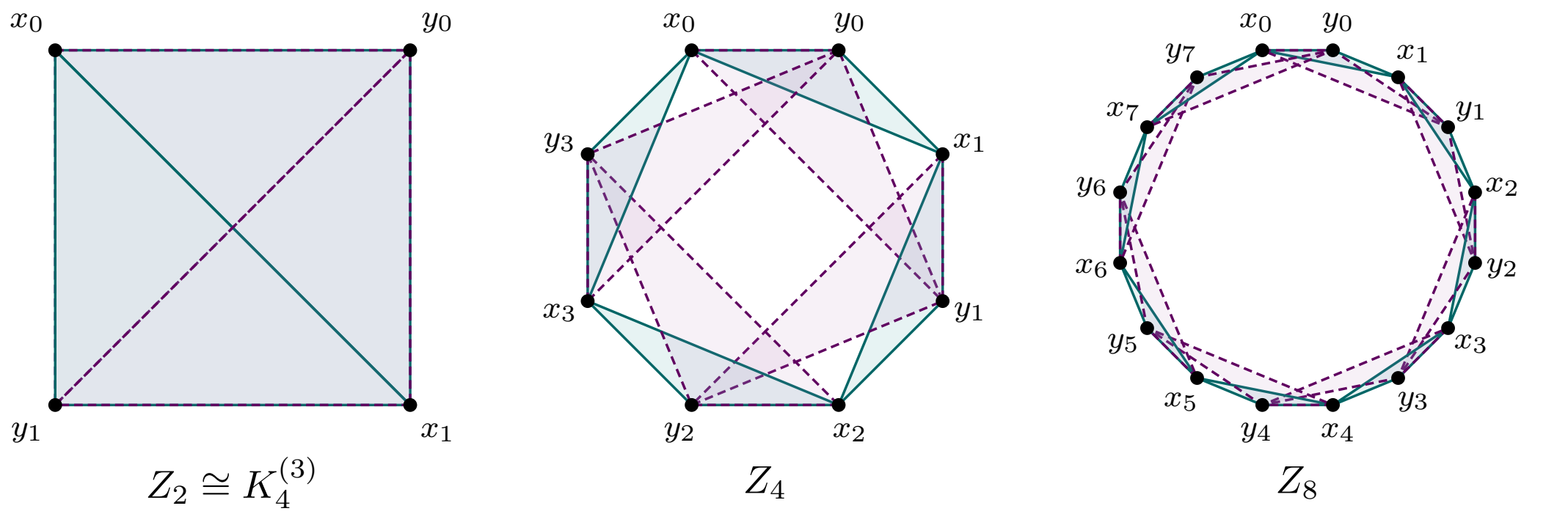

	In this paper, we determine the Tur\'an density of infinitely
	many members of another family of cyclic $3$-graphs.
	For an integer $\ell\ge2$, the \emph{zycle} $Z_\ell$ has
	vertex set and edge set
	$V(Z_\ell)=\{x_i,y_i:i\in\mathbb Z_\ell\}$, $
	E(Z_\ell)=\bigl\{\{x_i,y_i,x_{i+1}\},
	\{x_i,y_i,y_{i+1}\}:i\in\mathbb Z_\ell\bigr\}.$
	Here the $2\ell$ vertices are distinct and the indices are
	taken modulo $\ell$. Thus each designated pair
	$\{x_i,y_i\}$ forms an edge with each of the two vertices
	in the next pair. In particular, $Z_2\cong K_4^{(3)}$.
	The first examples relevant to our main result are shown
	in \Cref{fig:zycles}.
	
These hypergraphs were studied as \emph{rings} $R_\ell$
by DeBiasio and Jiang~\cite{DJ} in their work on the
co-degree threshold for the Fano plane.
They explicitly posed the problem of determining
the ordinary Tur\'an density or the codegree Tur\'an
density of each fixed zycle~\cite[Problem~1]{DJ}.

On the codegree side, Piga and Sch\"ulke~\cite{PS}
used zycles to prove that zero is an accumulation point
of the set of codegree Tur\'an densities of $3$-graphs.
In~\cite{WangCodegree}, we proved that the codegree
Tur\'an density of $Z_\ell$ is $1/3$ for every even
$\ell\ge8$ with $3\nmid\ell$, and at most $1/3$
for every even $\ell\ge8$.
These results give a partial answer to the codegree
part of DeBiasio and Jiang's problem, as well as
to a question of Piga and Sch\"ulke~\cite[Section~4]{PS}
concerning the sharpness of their upper bound of $1/2$.

In this paper, we address the ordinary Tur\'an density
part of DeBiasio and Jiang's problem.
For the family $\mathcal Z=\{Z_\ell:\ell\ge2\}$,
they proved that $\pi(\mathcal Z)=1/2$;
see~\cite[Theorem~4.4]{DJ}.
For individual members, they observed that
$\pi(Z_{pq})\le\pi(Z_p)$ for all integers $p\ge2$
and $q\ge1$, and that $\pi(Z_\ell)\ge5/9$
whenever $3\nmid\ell$;
see~\cite[Propositions~4.5--4.6]{DJ}.
Since $Z_2\cong K_4^{(3)}$, these bounds yield
\[
\frac59\le\pi(Z_\ell)\le\pi(K_4^{(3)})
\]
for every even integer $\ell\ge2$ with $3\nmid\ell$.

The lower bound is witnessed by the balanced
Tur\'an construction $T_n$, which has edge density
$5/9+o(1)$ and is $Z_\ell$-free whenever $3\nmid\ell$.
We recall this construction and verify the latter
assertion in \Cref{sec:density}.

Thus Tur\'an's tetrahedron conjecture,
$\pi(K_4^{(3)})=5/9$, would imply
$\pi(Z_\ell)=5/9$ for every even $\ell$
not divisible by three.
Our first result establishes this equality for every
such $\ell\ge4$, without assuming the conjecture,
and thereby gives a partial answer to
DeBiasio and Jiang's problem for ordinary
Tur\'an density.

	\begin{theorem}\label{thm:even}
		For every even integer $\ell\ge4$ with $3\nmid\ell$,
		$
		\pi(Z_\ell)=\frac59.
		$
		More generally, $\pi(Z_\ell)\le5/9$ for every even
		integer $\ell\ge4$.
	\end{theorem}
	
	We also determine the structure of almost extremal
	$Z_\ell$-free $3$-graphs in the range of equality.
	For two $3$-graphs on the same vertex set, their
	\emph{edit distance} is the size of the symmetric
	difference of their edge sets.
	
	\begin{theorem}\label{thm:stability}
		Fix an even integer $\ell\ge4$ with $3\nmid\ell$.
		For every $\varepsilon>0$, there exist $\delta>0$ and $n_0$
		such that every $Z_\ell$-free $3$-graph $G$ on $n\ge n_0$
		vertices with
		$
		e(G)\ge\left(\frac59-\delta\right)\binom n3
		$
		has edit distance at most $\varepsilon n^3$ from a copy of $T_n$.
	\end{theorem}
	
	\paragraph{Proof outline.}
	The upper bound uses the finite-family density theorem of
	Baber and Talbot, stated in \Cref{thm:BT-density}.
	We encode homomorphisms from zycles as closed walks in a
	digraph whose vertices are unordered pairs. Two short
	closed walks in a six-vertex $3$-graph yield the
	homomorphisms needed for every even length at least four.
	The deduction from the finite-family theorem is elementary;
	that theorem itself uses flag algebras.
	For stability, we combine the same homomorphisms with the
	hypergraph removal lemma and the Baber--Talbot stability
	theorem, using \Cref{lem:removal}  and \Cref{thm:BT-stability}.
	
	In \Cref{sec:preliminaries}, we record the homomorphism and
	blow-up tools. We prove \Cref{thm:even,thm:stability} in
	\Cref{sec:even}, and discuss the remaining even lengths
	in \Cref{sec:conclusion}.
	
	\section{Homomorphisms and blow-ups}\label{sec:preliminaries}
	
	A \emph{homomorphism} from an $r$-graph $F$ to an $r$-graph $G$
	is a map $\varphi:V(F)\to V(G)$ that sends each edge of $F$
	to an edge of $G$. In particular, the $r$ images of an edge
	are distinct. We write $F\to G$ when such a map exists.
	For an $r$-graph $G$ and a positive integer $t$, the
	\emph{blow-up} $G(t)$ is obtained by replacing each vertex
	of $G$ with an independent class of size $t$ and each edge
	with the complete $r$-partite $r$-graph on its corresponding
	classes. If $F\to G$, then
	$F\subseteq G(t)$ for $t\ge |V(F)|$.
	
	We use the following blow-up invariance theorem of Brown and
	Simonovits~\cite{BS}, in the finite-family form stated
	in~\cite[Theorem~7]{BT}.
	
	\begin{theorem}\label{thm:blowup}
		Let $\mathcal F=\{F_1,\ldots,F_s\}$ be a finite nonempty
		family of $r$-graphs, and let $t_1,\ldots,t_s$ be positive
		integers. Then
		\[
		\pi\bigl(\{F_1(t_1),\ldots,F_s(t_s)\}\bigr)=\pi(\mathcal F).
		\]
	\end{theorem}
	
	This allows us to transfer a density upper bound from a
	finite family to a single forbidden hypergraph.
	
	\begin{lemma}\label[lemma]{lem:transfer}
		Let $F$ be a $3$-graph and $\mathcal F$ a finite nonempty
		family of $3$-graphs. If $F\to J$ for every $J\in\mathcal F$,
		then $\pi(F)\le\pi(\mathcal F)$.
	\end{lemma}
	
	\begin{proof}
		Put $t=|V(F)|$. Since $F\subseteq J(t)$ for every
		$J\in\mathcal F$, every $F$-free $3$-graph is
		$\{J(t):J\in\mathcal F\}$-free. Apply \Cref{thm:blowup}.
	\end{proof}
	
	For zycles, the required homomorphisms can be expressed
	particularly simply. Given a pair $P\in\binom{V(G)}2$, let
	$
	N_G(P)=\{v\in V(G)\setminus P:P\cup\{v\}\in E(G)\}.
	$
	The \emph{pair digraph} $D(G)$ has vertex set
	$\binom{V(G)}2$, with
	\[
	P\to Q\quad\Longleftrightarrow\quad Q\subseteq N_G(P).
	\]
	Adjacent pair states are disjoint. All walks below are
	directed and may repeat vertices and arcs.
	
	\begin{lemma}\label[lemma]{lem:walk}
		For every $3$-graph $G$ and every integer $\ell\ge2$,
		$Z_\ell\to G$ if and only if $D(G)$ contains a closed walk
		of length $\ell$.
	\end{lemma}
	
	\begin{proof}
		Given a homomorphism $\varphi:Z_\ell\to G$, set
		$P_i=\{\varphi(x_i),\varphi(y_i)\}$. Each $P_i$ has two
		vertices, and the two edges at position $i$ give
		$P_{i+1}\subseteq N_G(P_i)$. Thus
		$P_0\to P_1\to\cdots\to P_{\ell-1}\to P_0$ is a closed walk.
		Conversely, assign $x_i,y_i$ bijectively to the two vertices
		of $P_i$ along such a walk. Every edge of $Z_\ell$ maps to
		an edge of $G$, so the resulting map is a homomorphism.
	\end{proof}
	
	\section{Proofs of the main results}\label{sec:even}
	
	\subsection{Tur\'an density}\label{sec:density}
	
	Let $H$ be the $3$-graph on $[6]$ with
	\begin{equation}\label{eq:H}
		E(H)=\{123,124,134,125,135,245,345,126,236,146,156,456\},
	\end{equation}
	where $ijk$ denotes the triple $\{i,j,k\}$.
	This is the graph $H_1$ in~\cite[Theorem~25]{BT}.
	We use the following density theorem of Baber and
	Talbot~\cite[Theorem~25]{BT}.
	
	\begin{theorem}\label{thm:BT-density}
		For the $3$-graph $H$ defined in \eqref{eq:H},
	$
		\pi\bigl(\{K_4^{(3)},H\}\bigr)=\frac59.
	$
	\end{theorem}
	
	We first construct homomorphisms from each even zycle to
	both members of this family. We abbreviate the pair
	$\{i,j\}$ by $ij$.
	
	\begin{lemma}\label[lemma]{lem:even-maps}
		For every even integer $\ell\ge4$,
		$Z_\ell\to K_4^{(3)}$ and $Z_\ell\to H$.
	\end{lemma}
	
	\begin{proof}
		In $D(K_4^{(3)})$, the pairs $12$ and $34$ form a directed
		two-cycle. Repeating it gives a closed walk of every even
		length, and hence $Z_\ell\to K_4^{(3)}$ by \Cref{lem:walk}.
		
		In $D(H)$, there are closed walks
		\begin{equation}\label{eq:short-walks}
			13\to45\to26\to13,
			\qquad
			13\to25\to14\to26\to13.
		\end{equation}
		Indeed, the successive arcs in the first walk use the edge
		pairs $(134,135)$, $(245,456)$ and $(126,236)$.
		The second uses $(123,135)$, $(125,245)$, $(124,146)$ and
		$(126,236)$. All these edges occur in \eqref{eq:H}.
		
		Both walks start at $13$, so they can be concatenated.
		Every even integer $\ell\ge4$ has the form $4a+6b$ with
		$a,b\ge0$: if $4\mid\ell$, take $(a,b)=(\ell/4,0)$;
		otherwise take $(a,b)=((\ell-6)/4,1)$.
		Concatenating $a$ copies of the four-step walk and $2b$
		copies of the three-step walk gives a closed walk of
		length $\ell$. Apply \Cref{lem:walk} again.
	\end{proof}
	
	We now turn to the matching lower bound. We use the
	classical Tur\'an construction; see, for example,~\cite{TetrahedronL2}.
	For pairwise disjoint sets $V_1,V_2,V_3$, let
	$T[V_1,V_2,V_3]$ be the $3$-graph on
	$V_1\cup V_2\cup V_3$ whose edges are all triples that
	either meet each part in one vertex, or have two vertices
	in $V_i$ and one in $V_{i+1}$ for some $i\in[3]$,
	where subscripts are taken modulo three.
	When the partition has $n$ vertices and its part sizes
	differ by at most one, we denote the resulting $3$-graph
	by $T_n$. Its number of edges is
	\begin{equation}\label{eq:turan-count}
		e(T_n)=|V_1||V_2||V_3|
		+\sum_{i=1}^3\binom{|V_i|}{2}|V_{i+1}|
		=\left(\frac59+o(1)\right)\binom n3.
	\end{equation}
	The following pair-type argument shows why divisibility
	by three enters the statement.
	
	\begin{lemma}\label[lemma]{lem:turan}
		Let $\ell\ge2$ with $3\nmid\ell$. For any pairwise disjoint
		sets $V_1,V_2,V_3$, there is no homomorphism from $Z_\ell$
		to $T[V_1,V_2,V_3]$. Consequently, $\pi(Z_\ell)\ge5/9$.
	\end{lemma}
	
	\begin{proof}
		Write $T=T[V_1,V_2,V_3]$ and consider a closed walk in $D(T)$.
		Call a pair of type $ii$ if both vertices lie in $V_i$,
		and of type $ij$ if it meets $V_i$ and $V_j$ for $i\ne j$.
		A pair contained in one part can only follow the transitions
		\[
		11\to22\to33\to11.
		\]
		Thus a walk cannot leave the collection of same-part pairs
		once it enters it. Among pairs meeting two parts, the only
		transitions that remain in that collection are
		\[
		12\to13\to23\to12.
		\]
		Indeed, the neighbourhoods of pairs of types $12$, $13$
		and $23$ are contained in $V_1\cup V_3$, $V_2\cup V_3$
		and $V_1\cup V_2$, respectively.
		A closed walk must stay in one of these two collections,
		and in either case its length is divisible by three.
		The first assertion follows from \Cref{lem:walk}.
		For balanced parts, \eqref{eq:turan-count} gives the
		claimed density lower bound.
	\end{proof}
	
	\begin{proof}[Proof of \Cref{thm:even}]
		For even $\ell\ge4$, \Cref{lem:transfer,lem:even-maps}
		and \Cref{thm:BT-density} give
		\[
		\pi(Z_\ell)\le\pi\bigl(\{K_4^{(3)},H\}\bigr)=\frac59.
		\]
		When $3\nmid\ell$, \Cref{lem:turan} gives equality.
	\end{proof}
	
	\subsection{Stability}
	
	To transfer stability, we first show that forbidding a fixed
	blow-up leaves only a negligible number of copies of its
	base graph. We use the following complete multipartite
	hypergraph theorem of Erd\H{o}s~\cite{Erdos}.
	
	\begin{theorem}\label{thm:complete-partite}
		For all integers $r\ge2$, $t\ge1$ and every $\eta>0$,
		there exists $n_0$ such that every $r$-graph on $n\ge n_0$
		vertices with at least $\eta n^r$ edges contains a complete
		$r$-partite $r$-graph with $t$ vertices in each part.
	\end{theorem}
	
	Applying this theorem to an auxiliary hypergraph whose
	edges are copies of a fixed $3$-graph gives the following
	counting consequence.
	
	\begin{lemma}\label[lemma]{lem:few-copies}
		Fix a $3$-graph $F$ without isolated vertices and a
		positive integer $t$. Every $F(t)$-free $3$-graph on $n$
		vertices has $o(n^{|V(F)|})$ labelled copies of $F$.
	\end{lemma}
	
	\begin{proof}
		Write $V(F)=[f]$. Otherwise, along a sequence of orders
		tending to infinity, there are at least $c n^f$ labelled
		copies for some fixed $c>0$. Randomly partition the vertices
		into $f$ labelled parts. Some partition retains at least
		$c n^f/f^f$ copies with vertex $i$ in part $i$.
		Form an $f$-partite $f$-uniform hypergraph whose edges are
		these copies. It has at least $c n^f/f^f$ edges, so
		\Cref{thm:complete-partite} gives a complete $f$-partite
		$f$-graph with $t$ vertices in each part. Its parts can
		be aligned with the original partition: any two vertices
		in different parts of this complete $f$-graph occur in a
		common edge and hence lie in different original parts.
		Thus we obtain $t$ vertices in each original part with
		every transversal an edge.
		For each edge $ijk$ of $F$, every triple on the chosen
		vertices in parts $i,j,k$ extends to a transversal copy
		of $F$. Hence all these triples are present, giving
		$F(t)$, a contradiction.
	\end{proof}
	
	The hypergraph removal lemma turns these small copy counts
	into a small set of edges whose deletion removes every
	forbidden copy. We use the following form; see
	R\"odl and Schacht~\cite[Theorem~2]{RS}.
	
	\begin{lemma}\label[lemma]{lem:removal}
		For every fixed $3$-graph $F$ and every $\varepsilon>0$,
		there exist $\delta>0$ and $n_0$ such that the following holds.
		If a $3$-graph $G$ on $n\ge n_0$ vertices contains at most
		$\delta n^{|V(F)|}$ labelled copies of $F$, then deleting
		at most $\varepsilon n^3$ edges makes $G$ $F$-free.
	\end{lemma}
	
	After these deletions, we can apply the stability theorem
	of Baber and Talbot~\cite[Theorem~26]{BT} for the same
	two-member family as in \Cref{thm:BT-density}.
	
	\begin{theorem}\label{thm:BT-stability}
		Let $(G_n)$ be a sequence of $\{K_4^{(3)},H\}$-free
		$3$-graphs with $|V(G_n)|=n$ and
		$e(G_n)=(5/9+o(1))\binom n3$.
		Then the edit distance from $G_n$ to some copy of $T_n$
		is $o(n^3)$.
	\end{theorem}
	
	We now combine these results with the homomorphisms
	constructed in \Cref{lem:even-maps}.
	
	\begin{proof}[Proof of \Cref{thm:stability}]
		It suffices to prove the sequential formulation.
		Let $(G_n)$ be a sequence of $Z_\ell$-free $3$-graphs with
		$e(G_n)=(5/9+o(1))\binom n3$.
		Put $t=2\ell$. By \Cref{lem:even-maps},
		\[
		Z_\ell\subseteq K_4^{(3)}(t),\qquad Z_\ell\subseteq H(t).
		\]
		Thus $G_n$ contains neither blow-up, and
		\Cref{lem:few-copies} implies that it has $o(n^4)$ copies
		of $K_4^{(3)}$ and $o(n^6)$ copies of $H$.
		
		Apply \Cref{lem:removal} to these two fixed graphs. Deleting $o(n^3)$ edges produces a
		$\{K_4^{(3)},H\}$-free subgraph $G_n'$ with
		$e(G_n')=(5/9+o(1))\binom n3$. \Cref{thm:BT-stability} gives a copy of $T_n$ at edit distance $o(n^3)$
		from $G_n'$, and hence also from $G_n$.
		
		If the stated $\varepsilon$--$\delta$ conclusion failed,
		there would be some $\varepsilon>0$ and a sequence of
		counterexamples of increasing orders with density at least
		$5/9-1/k$ at the $k$th term. By \Cref{thm:even}, their
		densities tend to $5/9$, contradicting the sequential
		conclusion just proved.
	\end{proof}
	
\section{Concluding remarks}\label{sec:conclusion}

\Cref{thm:even} determines $\pi(Z_\ell)$ for every even
$\ell\ge4$ not divisible by three. In particular,
\[
\pi(Z_4)=\pi(Z_8)=\frac59.
\]
The case $\ell=2$ remains outside the scope of our argument.
Indeed, $Z_2\cong K_4^{(3)}$ admits no homomorphism to $H$:
every such map would be injective, whereas the $3$-graph
$H$ in \eqref{eq:H} is $K_4^{(3)}$-free.
Thus the two-member family used here does not resolve
the tetrahedron problem.

Among even lengths at least four, the remaining cases are
precisely the multiples of six.
For these lengths, the classical Tur\'an construction
contains $Z_\ell$ for sufficiently large $n$, by placing
the designated pairs successively in $V_1,V_2,V_3$.
DeBiasio and Jiang~\cite[Theorem~4.4 and Proposition~4.5]{DJpre}
proved that $\pi(\mathcal Z)=1/2$ and
$\pi(Z_{pq})\le\pi(Z_p)$ for integers $p\ge2$ and $q\ge1$.
Together with their flag algebra bound
in~\cite[Proposition~4.10(iv)]{DJpre}, these results give
\[
\frac12\le\pi(Z_{6m})\le\pi(Z_6)<0.512303
\qquad\text{for every integer }m\ge1.
\]
This makes $Z_6$ the central remaining case: proving that
its density is $1/2$ would settle all multiples of six.
We conjecture that this is indeed the case.

\begin{conjecture}\label{conj:z6}
	$\pi(Z_6)=1/2$.
\end{conjecture}

By the preceding inequalities, Conjecture~\ref{conj:z6} would imply
$\pi(Z_\ell)=1/2$ whenever $6\mid\ell$.
Together with \Cref{thm:even}, this would determine
$\pi(Z_\ell)$ for every even $\ell\ge4$, leaving only
$Z_2$ among the even zycles.

Finally, our stability theorem describes the structure
of almost extremal $Z_\ell$-free $3$-graphs when
$\ell\ge4$ is even and $3\nmid\ell$.
It is natural to ask whether the balanced Tur\'an
construction is also exactly extremal for all
sufficiently large orders.

\begin{problem}
	For each even $\ell\ge4$ with $3\nmid\ell$, does
	$\ex(n,Z_\ell)=e(T_n)$ hold for all sufficiently large $n$?
\end{problem}


\begin{thebibliography}{99}
		
		\bibitem{BT}
		R.~Baber and J.~Talbot,
		New Tur\'an densities for $3$-graphs,
		\emph{Electron. J. Combin.} \textbf{19} (2012), no.~2, P22.
		\url{https://doi.org/10.37236/2360}.
		
		\bibitem{TetrahedronL2}
		L.~Bodn\'ar, W.~Chen, J.~Deng, J.~Hou, X.~Liu,
		J.~Song, J.~Yang and Y.~Zhang,
		Tetrahedron conjecture in the $\ell_2$-norm,
		preprint (2025). \arXiv{2511.12506}.
		
		\bibitem{BLLP}
		L.~Bodn\'ar, J.~Le\'on, X.~Liu and O.~Pikhurko,
		The Tur\'an density of short tight cycles,
		preprint (2025). \arXiv{2506.03223}.
		
		\bibitem{BS}
		W.~G.~Brown and M.~Simonovits,
		Digraph extremal problems, hypergraph extremal functions,
		and the densities of graph structures,
		\emph{Discrete Math.} \textbf{48} (1984), 147--162.
		
		\bibitem{dCF}
		D.~de Caen and Z.~F\"uredi,
		The maximum size of $3$-uniform hypergraphs not containing a Fano plane,
		\emph{J. Combin. Theory Ser. B} \textbf{78} (2000), no.~2, 274--276.
		\url{https://doi.org/10.1006/jctb.1999.1938}.
		
		\bibitem{DJ}
		L.~DeBiasio and T.~Jiang,
		On the co-degree threshold for the Fano plane,
		\emph{European J. Combin.} \textbf{36} (2014), 151--158.
		\url{https://doi.org/10.1016/j.ejc.2013.05.030}.
		
		\bibitem{DJpre}
		L.~DeBiasio and T.~Jiang,
		On the co-degree threshold for the Fano plane,
		expanded preprint, version~2 (2013).
		\arXiv{1204.1927v2}.
		
		\bibitem{Erdos}
		P.~Erd\H{o}s,
		On some extremal problems on $r$-graphs,
		\emph{Discrete Math.} \textbf{1} (1971), 1--6.
		\url{https://doi.org/10.1016/0012-365X(71)90002-1}.
		
		\bibitem{FF}
		P.~Frankl and Z.~F\"uredi,
		A new generalization of the Erd\H{o}s--Ko--Rado theorem,
		\emph{Combinatorica} \textbf{3} (1983), 341--349.
		\url{https://doi.org/10.1007/BF02579190}.
		
		\bibitem{FPS}
		Z.~F\"uredi, O.~Pikhurko and M.~Simonovits,
		The Tur\'an density of the hypergraph $\{abc,ade,bde,cde\}$,
		\emph{Electron. J. Combin.} \textbf{10} (2003), no.~1, R18.
		\url{https://doi.org/10.37236/1711}.
		
		\bibitem{KLP}
		N.~Kam\v{c}ev, S.~Letzter and A.~Pokrovskiy,
		The Tur\'an density of tight cycles in three-uniform hypergraphs,
		\emph{Int. Math. Res. Not. IMRN} \textbf{2024} (2024),
		no.~6, 4804--4841.
		
		\bibitem{KeevashSurvey}
		P.~Keevash,
		Hypergraph Tur\'an problems,
		in \emph{Surveys in Combinatorics 2011},
		London Math. Soc. Lecture Note Ser., vol.~392,
		Cambridge University Press, Cambridge, 2011, 83--140.
		\url{https://doi.org/10.1017/CBO9781139004114.004}.
		
		\bibitem{MubayiExpansion}
		D.~Mubayi,
		A hypergraph extension of Tur\'an's theorem,
		\emph{J. Combin. Theory Ser. B} \textbf{96} (2006),
		no.~1, 122--134.
		\url{https://doi.org/10.1016/j.jctb.2005.06.013}.
		
		\bibitem{PS}
		S.~Piga and B.~Sch\"ulke,
		Hypergraphs with arbitrarily small codegree Tur\'an density,
		\emph{Bull. Lond. Math. Soc.} \textbf{58} (2026), no.~4, e70348.
		\url{https://doi.org/10.1112/blms.70348}.
		
		\bibitem{PikhurkoExpansion}
		O.~Pikhurko,
		Exact computation of the hypergraph Tur\'an function
		for expanded complete $2$-graphs,
		\emph{J. Combin. Theory Ser. B} \textbf{103} (2013),
		no.~2, 220--225.
		\url{https://doi.org/10.1016/j.jctb.2012.09.005}.
		
		\bibitem{RS}
		V.~R\"odl and M.~Schacht,
		Generalizations of the removal lemma,
		\emph{Combinatorica} \textbf{29} (2009), 467--501.
		
		\bibitem{WangCodegree}
		S.~Wang,
		Exact codegree Tur\'an densities of even zycles,
		preprint (2026). \arXiv{2611.12582}.
		
	\end{thebibliography}
\end{document}